\documentclass[a4paper,oneside,reqno,11pt]{amsart}
\usepackage[utf8]{inputenc}
\usepackage[english]{babel}
\usepackage[margin=3.6cm]{geometry}
\usepackage{graphicx}
\usepackage{amsmath,amssymb,amsthm}
\usepackage{enumerate}
\usepackage{microtype}
\newtheorem{theorem}{Theorem}[section]
\newtheorem{lemma}{Lemma}[section]
\newcommand{\rr}{\mathbf{R}}
\newcommand{\mop}[1]{\mathop{\mathrm{#1}}}
\newcommand{\ocurl}{\mop{curl}}
\newcommand{\odiv}{\mop{div}}
\newcommand{\otr}{\mop{tr}}
\newcommand{\Idx}[2]{\int_{#1} #2 \,dx}

\newcommand{\IdX}[3]{\int_{#1} #2 \,d#3}
\newcommand{\Inv}[1]{\frac{1}{#1}}
\newcommand{\inv}{^{-1}}
\newcommand{\ddt}[2]{\frac{\partial #1}{\partial #2}}
\newcommand{\lrp}[1]{\left(#1\right)}
\newcommand{\lran}[1]{\langle#1\rangle}
\newcommand{\e}{{\varepsilon}}
\newcommand{\mat}[1]{\begin{pmatrix}#1\end{pmatrix}}
\newcommand{\figthreetex}[6]{ % three figures side by side
  \begin{figure}[#6]
  \begin{center}
  \begin{minipage}[c]{.31\textwidth}
  \centering\includegraphics[width=.95\textwidth]{#3}
  \end{minipage}
  \begin{minipage}[c]{.01\textwidth}
  $\,$
  \end{minipage}
  \begin{minipage}[c]{.31\textwidth}
  \centering\includegraphics[width=.95\textwidth]{#4}
  \end{minipage}
  \begin{minipage}[c]{.01\textwidth}
  $\,$
  \end{minipage}
  \begin{minipage}[c]{.31\textwidth}
  \centering\includegraphics[width=.95\textwidth]{#5}
  \end{minipage}
  \begin{minipage}[b]{.31\textwidth}
  \vspace{.3cm}\centering(a)
  \end{minipage}
  \begin{minipage}[c]{.01\textwidth}
  $\,$
  \end{minipage}
  \begin{minipage}[b]{.31\textwidth}
  \vspace{.3cm}\centering(b)
  \end{minipage}
  \begin{minipage}[c]{.01\textwidth}
  $\,$
  \end{minipage}
  \begin{minipage}[b]{.31\textwidth}
  \vspace{.3cm}\centering(c)
  \end{minipage}
  \caption{#2}
  \label{#1}
  \end{center}
  \end{figure}}
\newcommand{\figgtex}[3]{
  \begin{figure}[hb]
  \begin{center}
  \includegraphics[width=.5\textwidth]{#3}
  \end{center}\caption{#2}\label{#1}
  \end{figure}}

\author{Dag Lukkassen$^{\star\dagger}$}
\address{$^{\star}$Narvik University College, P.O. Box 385, N-8505 Narvik, Norway}
\address{$^{\dagger}$Norut Narvik, P.O. Box 250, N-8504 Narvik, Norway}
\author{Annette Meidell$^{\star\dagger}$}
\author{Klas Pettersson$^{\star}$}
\email{klapet@hin.no}

\date{July 16, 2015}

\title{A local-effective relation in planar linear elasticity}

\begin{document}

\begin{abstract}
We give an example of a relation between local and 
  effective properties for elastic structures, up to 
  geometric constants.
The model considered is a periodic structure with isotropic and homogeneous
 local elasticity tensor in planar linear elasticity.
The corresponding physical model is a flat two dimensional body with holes as, for example, a perforated plate.
\end{abstract}

\maketitle

\section{Introduction}

In this paper we consider the effective elastic properties of a periodic structure in planar linear elasticity.
The structure is assumed to be two dimensional with stress-free holes and locally isotropic.
We show that the relation between local and effective elastic properties is elementary under natural hypotheses.
The result is an extension of previous work by Vigdergauz \cite{vigdergauz}, which is closely related to the CLM theorem \cite{milton1}, and gives an example of a so-called link in the general theory described in \cite{miltonbook}.

In \cite[§44]{muskhelishvili} Muskhelishvili has given a characterization of the dependence of the state of stress on the elastic moduli for homogeneous isotropic bodies in the planar theory of linear elasticity. He has shown that 
for multiply connected domains with external forces applied to each connected part of the boundary which are separately zero, the state of stress does not depend on the elastic moduli. 
If on the other hand there are holes with nonzero external forces,
the state of stress do depend on the value of $\lambda/\mu$, where $\lambda$ and $\mu$ are the Lam{\'e} moduli.
Here it seems that attributions could be made to Levy \cite{levy1898legitimite} and Michell \cite{michell}.

In \cite{olverelast} Olver has shown that every linear planar anisotropic elastic material is equivalent, under a linear change of coordinates, to an orthotropic material. In particular, he has shown that two isotropic Lagrangians determine the same orthotropic Lagrangian if and only if their {Lam\'e} moduli are proportional $\lambda/\mu = \tilde \lambda/ \tilde \mu$.

In \cite{milton1} Cherkaev, Lurie, and Milton have investigated the dependence of effective elastic properties of composite materials on the value of the Poisson ratio $\nu = (1 + \mu/\lambda)\inv$.
They have shown that a uniform shift in the Poisson ratio gives the same shift in the effective elastic properties of the material under the assumption of effective isotropy.

In \cite{vigdergauz} Vigdergauz has extended the result of \cite{milton1} to planar elastic structures with effective square symmetry. In particular,
he has shown that 
\begin{align*}
\frac{1}{K^*} & = \frac{1}{K} + A_1\Big(\frac{1}{K} + \frac{1}{\mu}\Big) \\
\frac{1}{\mu^*} & = \frac{1}{\mu} + A_2\Big(\frac{1}{K} + \frac{1}{\mu}\Big) \\
\frac{1}{\mu^*_{45}} & = \frac{1}{\mu} + A_3\Big(\frac{1}{K} + \frac{1}{\mu}\Big),
\end{align*}
where $A_i$ are geometric constants, $K$ denotes the local bulk modulus, $K^*$ is the effective bulk modulus, and $\mu^*$, $\mu^*_{45}$ are the effective shear moduli.

We give an extension of Vigdergauz's result 
to the case when the effective elastic tensor is allowed to be
anisotropic.
This leads to six geometric constants instead of the three for effective square symmetry.
We show that the inverse $s^*_{ijkl}$ of the effective elasticity tensor $c^*_{ijkl}$ satisfies the equation (Theorem \ref{tm:main})
\begin{align}\label{eq:result}
s^*_{ijkl} & = s_{1111} d_{ijkl} + s_{1212} e_{ijkl},
\end{align}
where $d_{ijkl}$ are constants which do not depend on the local elasticity tensor $c_{ijkl}$,
and $e_{ijkl}$ are constants independent of both geometry and elastic moduli.
Like the CLM and Vigdergauz theorems, the relation \eqref{eq:result} separates the dependence of the effective elasticity tensor on the geometry from the dependence on the local elastic properties in the material.

The proof of \eqref{eq:result} we give below is based on the stress invariance of Muskhelishvili.
This approach has been illustrated in the note \cite{lukkassen} for a class of effectively square symmetric structures.

\section{Statement of result}

In a bounded Lipschitz domain $\Omega$ in $\rr^2$ we consider the total elastic energy corresponding to a displacement $u$ given by
\begin{align}\label{eq:I}
I(u) & = \int_{\Omega} W(\nabla u(x)) \,dx,
\end{align}
where $W$ is a symmetric quadratic stored energy function of the material:
\begin{align*}
W(\nabla u) & = \sum_{i,j,k,l} a_{ijkl}\frac{\partial u_i}{\partial x_j}\frac{\partial u_k}{\partial x_l}.
\end{align*}
The variational moduli $a_{ijkl}$ are assumed to satisfy $a_{ijkl} = a_{klij}$.
By the symmetry of $W$ there are elastic moduli $c_{ijkl}$ such that 
\begin{align*}
W(\nabla u) & = \frac{1}{2}\sum_{i,j,k,l} c_{ijkl} \varepsilon_{ij}(u) \varepsilon_{kl}(u),
\end{align*}
and
\begin{align}\label{eq:tricl}
c_{ijkl} & = c_{jikl} = c_{klij},
\end{align}
where $\varepsilon$ denotes the symmetric part of the gradient.
The elastic moduli are assumed to be bounded and satisfy the Legendre-Hadamard strong ellipticity condition: There are positive constants $\kappa_i$ such that
\begin{align*}
\kappa_1 |\xi|^2 \le \sum_{i,j,k,l} c_{ijkl}(x) \xi_{ij} \xi_{kl} \le \kappa_2 |\xi|^2,
\end{align*}
for almost every $x$ in the domain and every symmetric $\xi$.
The Euler-Lagrange equations corresponding to \eqref{eq:I} are 
\begin{align*}
\sum_{j,k,l}a_{ijkl}\frac{\partial^2 u_k}{\partial x_j \partial x_l} & = 0,
\end{align*}
which in terms of the elastic moduli may be written
\begin{align*}
\mop{div}\sigma(u) & = 0,
\end{align*}
where
\begin{align*}
\sigma_{ij}(u) = \sum_{k,l}c_{ijkl}\varepsilon_{kl}(u)
\end{align*}
is the Hooke law.

The material is homogeneous if $c_{ijkl}$ are constant on $\Omega$.
The material is isotropic if the stored energy function has the form
\begin{align*}
W(\nabla u) & = \frac{1}{2}\sum_{i,j,k,l}c_{ijkl}\varepsilon_{ij}(u)\varepsilon_{kl}(u)
= \mu |\varepsilon(u)|^2 + \frac{\lambda}{2}(\mop{tr}\varepsilon(u))^2,
\end{align*}
where $\mu$ and $\lambda$ are the Lam\'e moduli; $\mop{tr}$ denotes the trace.
In addition to \eqref{eq:tricl} the components of an isotropic tensor thus satisfy
\begin{align*}
c_{1111} & = c_{2222} = c_{1122} + 2c_{1212}, &
c_{1112} & = c_{2212} = 0.
\end{align*}

Let $\Omega_h$ be a periodically perforated domain constructed as follows.
Let $Q$ be a bounded and connected domain contained in some translation $Y$ of the periodicity cell $(0,l_1) \times (0,l_2)$.
We assume that the periodically extended domain $\widetilde{Q} = \mop{Int} \mop{Cl} ( Q + l_1 \mathbf{Z} \times l_2 \mathbf{Z} )$
is connected and Lipschitz.
For $h > 0$, we set $\Omega_h$ to be interior of the union of all 
 translates $\mop{Cl}(Q + (l_1 k_1, l_2 k_2))$ contained in $\Omega$, $k \in \mathbf{Z}^2$.

In the periodically perforated domain $\Omega_h$ we consider the total elastic energy 
\begin{align}\label{eq:Ive}
I^h(u) & = \frac{1}{2}\int_{\Omega_h} \sum_{i,j,k,l} c_{ijkl}(\frac{x}{h})\varepsilon_{ij}(u)\varepsilon_{kl}(u) \,dx,
\end{align}
where the elastic moduli are assumed to be $Y$-periodic.

We say that a vector $u \in H^1(Q)^2$ is quasiperiodic if 
  $u - \xi x \in H^1_{\mop{per}}(Q)^2$ for some constant matrix
  $\xi$, where $H^1_{\mop{per}}(Q)$ denotes the closure in
  $H^1(Q)$ of the smooth periodic functions
  $C^\infty_{\mop{per}}(Y)$.
We call $\xi$ the quasiperiod of $u$.

In the asymptotic limit as $h$ tends to zero, the functional $I^h$ behaves as
the homogenized functional
\begin{align*}
I^*(u) & = \frac{1}{2}\int_{\Omega} \sum_{i,j,k,l} c^*_{ijkl}\varepsilon_{ij}(u)\varepsilon_{kl}(u) \,dx,
\end{align*}
where the effective elastic moduli $c^*_{ijkl}$ are homogeneous and defined for symmetric $\xi$ by
\begin{align}\label{eq:ceff}
c^*_{ijkl} \xi_{kl} & = \frac{1}{|Y|}\int_Q c_{ijkl}\varepsilon(u^\xi) \, dx,
\end{align}
where $u^\xi$ is quasiperiodic in $Q$ with quasiperiod $\xi$ and a  minimizer of the total cell energy
\begin{align*}
J(u) & = \frac{1}{2}\int_Q c_{ijkl}\varepsilon_{ij}(u)\varepsilon_{kl}(u) \,dx.
\end{align*}

The effective tensor defined by \eqref{eq:ceff} is by the strong ellipticity and periodicity, positive definite
on the set of symmetric matrices \cite{oleinik}.

The asymptotic behavior of the functional $I^h(u)$ may be illustrated by the Dirichlet
boundary value problem
\begin{align*}
\mop{div}\sigma(u) + f & = 0 \text{ in } \Omega_h,\\
u & = 0 \text{ on } \partial \Omega.
\end{align*}
The sequence of solutions $u_h \in H^1_0(\Omega_h, \partial \Omega)$ to the above problem
admits an asymptotic expansion of the form
\begin{align*}
u_0 + h \sum_i N^i\big(\frac{x}{h}\big) \frac{\partial u_0}{\partial x_i}
\end{align*}
in the sense that 
\begin{align*}
\big\| u_h - u_0 - h \sum_i N^i\big(\frac{x}{h}\big) \frac{\partial u_0}{\partial x_i} \big\|_{H^1(\Omega_h)} \to 0,
\end{align*}
as $h$ tends to zero, where $u_0$ is the minimizer of $I^*(u) + |Q||Y|\inv\int_\Omega f u dx$ on $H^1_0(\Omega)$, and $N^i$ are auxiliary periodic functions.
See \cite{jikov,oleinik}.

The first assertion (a) in the following theorem is the result of this paper.
For isotropic effective tensors it has been observed numerically in~\cite{day1992elastic}.
The second part (b) is a version of the CLM theorem \cite{milton1},
which is covered by the same argument of proof.

\begin{theorem}\label{tm:main}
Let $c_{ijkl}$ be a homogeneous and isotropic elasticity tensor in the perforated periodicity cell $Q$.
Let $c^*_{ijkl}$ be the effective elasticity tensor defined by \eqref{eq:ceff} with local tensor $c_{ijkl}$.
Let $s_{ijkl}$ and $s^*_{ijkl}$ be the inverses of $c_{ijkl}$ and $c^*_{ijkl}$, respectively.
Then there exist constants $d_{ijkl}$, which are independent of $c_{ijkl}$, such that
\begin{align}
s^*_{ijkl} & = s_{1111} \, d_{ijkl} + s_{1212} \, e_{ijkl},\tag{a}
\end{align}
where $e_{ijkl}$ satisfies $e_{ijkl} = e_{jikl} = e_{klij}$ and are defined
by $e_{1122} = -2e_{1212} = -2$, with the rest of the components equal to zero.

Let $c_{ijkl}^{(1)}$ and $c_{ijkl}^{(2)}$ be two smooth isotropic elasticity tensors on $Q$.
Let $c^{*(1)}_{ijkl}$ and $c^{*(2)}_{ijkl}$ be the effective elasticity tensors defined by 
\eqref{eq:ceff} with local tensors $c_{ijkl}^{(1)}$ and $c_{ijkl}^{(2)}$, respectively.
Let $s_{ijkl}^{(1)}$, $s_{ijkl}^{(2)}$, $s^{*(1)}_{ijkl}$, and $s^{*(2)}_{ijkl}$ be the inverses of
$c_{ijkl}^{(1)}$, $c_{ijkl}^{(2)}$, $c^{*(1)}_{ijkl}$, and $c^{*(2)}_{ijkl}$, respectively.
Suppose that
\begin{align*}
s_{1122}^{(1)} - s_{1122}^{(2)} = -2(s_{1212}^{(1)} - s_{1212}^{(2)}) = \mathrm{constant}.
\end{align*}
Then
\begin{align}
s^{*(1)}_{ijkl} - s^{*(2)}_{ijkl} = (s^{(1)}_{1212} - s^{(2)}_{1212})e_{ijkl} = s^{(1)}_{ijkl} - s^{(2)}_{ijkl}.\tag{b}
\end{align}
\end{theorem}

\vspace{.25cm}

A proof of Theorem~\ref{tm:main} is presented in Section~\ref{sec:theproof}.
In Section~\ref{sec:aux} we prove the lemmas stated in Section \ref{sec:theproof}.

\section{Proof of Theorem~\ref{tm:main}}\label{sec:theproof}

Theorem \ref{tm:main} is concerned with the effective elasticity tensor defined by \eqref{eq:ceff}.
In the following we will therefore only work in the periodicity cell $Q$ or a periodic extension of it.
An exception is Lemma \ref{lm:michell}, which is of a more general nature.

We choose curves connecting opposite sides of the perforated periodicity cell $Q$.
Let $\gamma_1$ be a smooth curve connecting two points in $\mop{Cl}Q$
with $l_1$ as the difference between the $x_1$ components, and with equal
$x_2$ components. We assume that $\gamma_1$ is separated from any possible
hole in the global structure, that is
$\mop{dist}(\gamma_1,\partial \widetilde{Q}) > 0$.
Moreover, we suppose that the $x_1$ component of the unit tangent $\tau$ to
$\gamma_1$ is always positive.
Such a curve $\gamma_1$ exists by the connectedness of $Q$ and
its periodic extension $\widetilde{Q}$.
Let $\gamma_2$ be an analogue curve in the $x_2$ direction, by interchanging
all the indices.

For scalar valued $\phi$ and vector valued $\varphi$, we denote 
\begin{align*}
  \mop{curl} \phi & = \lrp{ \ddt{\phi}{x_2}, -\ddt{\phi}{x_1} }, &
  \mop{curl} \varphi & = \ddt{\varphi_2}{x_1} - \ddt{\varphi_1}{x_2}.
\end{align*}

We choose the following curves in the extended domain.
Let $\gamma_1'$ be the translate of $\gamma_1$ such that the left endpoint of
  $\gamma_1'$ is the common point of $\gamma_1$ and $\gamma_2$.
Then let $\gamma_2'$ be the translate of $\gamma_2$ such that its lower endpoint
  is the left endpoint of $\gamma_1'$.
Let $Q'$ be the intersection of $\widetilde{Q}$ and the region of $\rr^2$ bounded by
  \[
  \Gamma := \gamma_1' \cup \gamma_2' \cup ((0, l_2) + \gamma_1' )\cup ((l_1, 0) + \gamma_2').
  \]
The domains and the curves are illustrated in Figures \ref{fig:structure} and \ref{fig:cells}(a)--(c).

When there is only one displacement field $u$ and one elasticity tensor $c_{ijkl}$ under consideration
we will write $\sigma$ in place of $\sigma(u)$ where $\sigma_{ij}(u) = \sum_{k,l}c_{ijkl}\varepsilon_{kl}(u)$.
Similarly, $\varepsilon = \varepsilon(u)$ and $\omega = \omega(u)$, where $\omega$ denotes the
antisymmetric part of the gradient: $\nabla = \varepsilon + \omega$.

\begin{lemma}\label{lm:qpspecial}
Let $c_{ijkl}$ be a constant isotropic elasticity tensor on $Q$, periodically extended.
Suppose that $u$ is quasiperiodic on $Q$ with symmetric quasiperiod $\xi$.
Assume that $\mop{div}\sigma = 0$ on $\gamma_1 \cup \gamma_2$.
Suppose that $c_{ijkl}$ and $u$ are smooth in some neighborhood of $\gamma_i$.
Then
\begin{align*}
\xi_{11} & = -  \frac{2s_{1212}}{|Y|} \int_Q  \sigma_{22} \,dx + \frac{s_{1111}}{l_1}  \int_{\gamma_1}  ( (\otr \sigma, 0) - x_2 \ocurl \otr \sigma ) \,dx, \\
\xi_{22} & = -  \frac{2s_{1212}}{|Y|} \int_Q  \sigma_{11} \,dx + \frac{s_{1111}}{l_2}  \int_{\gamma_2}  ( (0, \otr \sigma) + x_1 \ocurl \otr \sigma ) \,dx, \\
\xi_{12} & = \frac{2s_{1212}}{|Y|} \int_Q \sigma_{12} \,dx + \frac{s_{1111}}{2l_1} \int_{\gamma_1'}  ( (\otr \sigma, 0) + x_1 \ocurl \otr \sigma ) \,dx  \\
     & \qquad\qquad\qquad\qquad\,\, + \frac{s_{1111}}{2l_2} \int_{\gamma_2'}  ( (0,\mop{tr}\sigma) - x_2 \ocurl \mop{tr}\sigma ) \,dx,
\end{align*}
where $s_{ijkl}$ is the inverse of $c_{ijkl}$.
\end{lemma}

For pairs of elasticity tensors $c_{ijkl}^{(1)}$, $c_{ijkl}^{(2)}$, and the corresponding
inverses, stresses, etc, we will use the prefix $\Delta$ to denote the difference.
For example, $\Delta c_{ijkl} = c_{ijkl}^{(1)} - c_{ijkl}^{(2)}$.

\begin{lemma}\label{lm:qpspecial2}
Let $c_{ijkl}$ be an isotropic elasticity tensor on $Q$, periodically extended.
Let $u$ be quasiperiodic on $Q$ with symmetric quasiperiod $\xi$.
Suppose that $c_{ijkl}$ and $u$ are smooth in some neighborhood of $\gamma_i$.
  Assume that $\mop{div}\sigma = 0$ on $\gamma_1 \cup \gamma_2$.
Let $s_{ijkl}$ denote the inverse of $c_{ijkl}$.
Assume that $s_{ijkl}$ is shifted in such a way that
$\Delta s_{1122} = -2\Delta s_{1212} = \mathrm{constant}$.
Then
\begin{align*}
\Delta \xi  & = 
- \frac{2\Delta s_{1212}}{|Y|} \int_Q \mat{ \sigma_{22} & - \sigma_{12} \\ - \sigma_{12} & \sigma_{11} } \,dx.
\end{align*}
\end{lemma}

\figgtex{fig:structure}{An interior part of a periodically perforated structure.}{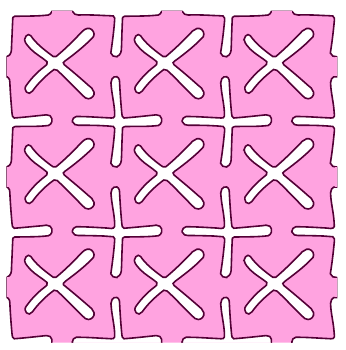}

\figthreetex{fig:cells}
  {(a) A periodicity cell $Q$ with curves $\gamma_1$ and $\gamma_2$.
   (b) Translations $\gamma_1'$ and $\gamma_2'$ of $\gamma_1$ and $\gamma_2$,
         respectively.
   (c) The domain $Q'$, with outer boundary $\Gamma$ corresponding to
         the curves in (a) and (b); $Q'$ has no cusps.}
  {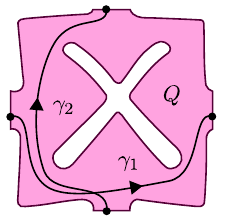}
  {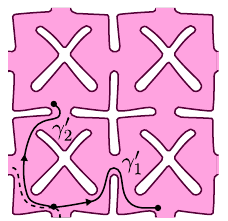}
  {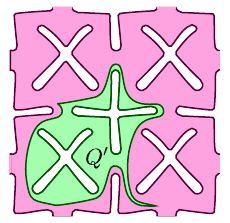}
  {hbp}

The following result is concerned with the dependence of the stress on the elasticity tensor for the Neumann problem. See \cite[§44]{muskhelishvili}.

\begin{lemma}\label{lm:michell}
Let the elasticity tensor $c_{ijkl}$ be smooth and isotropic on a bounded and connected Lipschitz domain $\Omega$ in $\rr^2$.
Then the stress $\sigma$ does not depend on the elasticity tensor in the pure
Neumann problem on $\Omega$
with $\lran{\sigma \nu, 1}_{\Gamma} = 0$ for all connected components $\Gamma$ of $\partial \Omega$,
\begin{enumerate}[(i)]
\item if $c_{ijkl}$ is constant, or
\item if $\Delta s_{1122} = -2\Delta s_{1212} = \mathrm{constant}$,
\end{enumerate}
where $s_{ijkl}$ denotes the inverse of $c_{ijkl}$.
In (i) the condition on the data is necessary if $\Delta(c_{ijkl}/c_{1111}) \neq 0$.
On the other hand, if $c_{ijkl}$ is constant and $\Delta(c_{ijkl}/c_{1111}) = 0$,
the condition on the data can be relaxed to just equilibrium, $\lran{\sigma \nu, 1}_{\partial \Omega} = 0$,
and the stress does not depend on the specific components of the elasticity tensor.
\end{lemma}

\begin{proof}[Proof of Theorem~\ref{tm:main}]
Since $c^*_{ijkl}$ defined by \eqref{eq:ceff} is positive definite,
we have that for every average stress there exist a unique stress field $\sigma$ that, by Lemma~\ref{lm:michell}, does not depend on the
isotropic local elasticity tensor $c_{ijkl}$.
This is because the boundary of any hole in the global problem is stress-free
and hence the average normal stress on the outer part $\Gamma_0$ of the boundary of $Q$ vanishes by the Green formula. 
The corresponding quasiperiodic displacement field $u$ is unique up to
a translation.
The symmetric quasiperiod $\xi$ of $u$ is therefore uniquely determined by $\sigma$ and $c_{ijkl}$.

Let $\lran{\sigma}$ denote the average stress
$|Y|\inv \int_Q (\sigma_{11},\sigma_{22},\sigma_{12}) \,dx$, which  
we will let vary over the canonical basis of $\rr^3$. 
With $\lran{\sigma} = (1,0,0)$, \eqref{eq:ceff}
gives
$
s^*_{1111} = ((c^*)\inv)_{1111} = \xi_{11}. 
$
Since $c_{ijkl}$ is isotropic, $\xi$ is symmetric, and $\mop{div}\sigma$ vanishes,
there exists by
Lemma~\ref{lm:qpspecial} a real constant $d_{1111}$, which
depends only on $\sigma$ and hence not on $c_{ijkl}$, such that
$s^*_{1111} = s_{1111}\, d_{1111}$.
Moreover, by the same equation, there exists constants $d_{1122}$ and $d_{1112}$ such that
$
s^*_{1122} = -2s_{1212} + s_{1111} \, d_{1122}
$
and
$
s^*_{1112} = s_{1111}\, d_{1112}.
$
With $\lran{\sigma} = (0,1,0)$, we similarly find constants $d_{2222}$ and
$d_{2212}$ such that
$s^*_{2222} = s_{1111}\, d_{2222}$ and
$s^*_{2212} = s_{1111}\, d_{2212}$.
Finally, with $\lran{\sigma} = (0,0,1)$, we find a constant $d_{1212}$ such that
$
s^*_{1212} = s_{1212} + s_{1111}\, d_{1212}.
$
Let $d$ have the entries $d_{ijkl}$ defined above such that it satisfies the symmetries 
$d_{ijkl} = d_{jikl} = d_{klij}$.
Let $e_{ijkl}$ be as in the statement of the theorem.
Then $s^* = s_{1111}\, d + s_{1212} \, e$,
where $d_{ijkl}$ and $e_{ijkl}$ do not depend on $c_{ijkl}$ by construction.

We argue in the same way for the case when $c_{ijkl}$ is changed in such a way that
  $\Delta s_{1122} = -2\Delta s_{1212} = \mathrm{constant}$.
First, we have $\Delta s_{1111} = 0$ by isotropy,
  and by the assumption of uniform shift of the not necessarily constant
  compliance tensor $s_{ijkl}$, the stress is invariant by (ii) in Lemma~\ref{lm:michell}.
By Lemma~\ref{lm:qpspecial2}, we have
  $\Delta s^*_{1111} = \Delta \xi_{11} = 0$,
  $\Delta s^*_{1122} = -2 \Delta s_{1212} = \Delta s_{1122}$, and
  $\Delta s^*_{1112} = 0$,
  for $\lran{\sigma} = (1,0,0)$.
For $\lran{\sigma} = (0,1,0)$, we find
  $\Delta s^*_{2222} = 0$ and $\Delta s^*_{2212} = 0$.
Finally, we find $\Delta s^*_{1212} = \Delta s_{1212}$ for $\lran{\sigma} = (0,0,1)$.
In conclusion, $\Delta s^* = \Delta s_{1212} \, e = \Delta s$.
\end{proof}

\section{Proofs of Lemmas 3.1--3}\label{sec:aux}

In this section we prove Lemma~\ref{lm:qpspecial}--\ref{lm:michell}.
We will use the following Ces\`{a}ro formula.
Let $\varsigma_{ij}$ denote $0$ if $i = j$, and $1$ otherwise.

\begin{lemma}\label{lm:quasiperiodlemma}
Let $c_{ijkl}$ be an elasticity tensor.
Let $u$ be quasiperiodic on $Q$ with quasiperiod $\xi$.
Suppose that $c_{ijkl}$ and $u$ are smooth in some neighborhood of $\gamma_i$.
Then
  \begin{align*}
  & \xi_{ij} + \omega_{ji}(\alpha) \\ 
  & \!\!\! = \Inv{l_i}
  \sum_{q,r,s,t} \Bigl( 
  \frac{1}{2}\!\!\int_{\gamma_i} s_{jrst} \sigma_{st} \,dx_r
  + \varsigma_{rj}  \!\int_{\gamma_i} \ddt{(s_{rqst}\sigma_{st})}{x_j} x_r \,dx_q
  - \varsigma_{rj}  \!\int_{\gamma_i} \ddt{(s_{jqst}\sigma_{st})}{x_r} x_r \,dx_q \Bigr), 
  \end{align*}
  where $s_{ijkl}$ is the inverse of $c_{ijkl}$, and $\alpha$ is the starting point of the integration along
  $\gamma_i$.
\end{lemma}
\begin{proof}
Since $u$ is quasiperiodic and smooth in a neighborhood of $\gamma_i$,
  and the endpoints $\alpha$ and $\beta$ of $\gamma_i$ have equal components except
  for the $i$ths, we have by the Newton-Leibniz formula,
\begin{align*}
\xi_{ij}l_i & = u_j(\beta) - u_j(\alpha)
 = \lrp{ \Idx{\gamma_i}{\nabla u} }_j
 = \lrp{ \Idx{\gamma_i}{\e} }_j
   + \lrp{ \Idx{\gamma_i}{\omega} }_j.
\end{align*}
By the linear relation between $\e$ and $\sigma$, and the nonsingularity of the elasticity tensor $c_{ijkl}$, we have
  \begin{align}\label{eq:auxqp}
  \lrp{ \Idx{\gamma_i}{\e} }_j
  & =
  \IdX{\gamma_i}{\sum_r \e_{jr}}{x_r}
  =
  \sum_{r,s,t} \IdX{\gamma_i}{s_{jrst}\sigma_{st}}{x_r}.
  \end{align}
By the smoothness of $u$, an integration by parts yields
  \begin{align*}
  \lrp{ \int_{\gamma_i} \omega \,dx }_j  
  & =
  l_i \omega_{ji}(\alpha) + \sum_{\substack{k,r \\ r\neq j}}\int_{\gamma_i} 
  \ddt{\omega_{rj}}{x_k}
  x_r
  \,dx_k,
  \end{align*}
  where we used that $\omega$ is periodic since $u$ is quasiperiodic.
After using that
  \begin{align*}
  \ddt{\omega_{rj}}{x_k}
  & =
  \ddt{\e_{rk}}{x_j} - \ddt{\e_{jk}}{x_r},
  \end{align*}
  one finds
  \begin{align*}
  \int_{\gamma_i} 
  \ddt{\omega_{rj}}{x_k}
  x_r
  \,dx_k 
  & =
  \sum_{s,t}\Bigl( 
  \int_{\gamma_i}  
  \ddt{(s_{rkst}\sigma_{st})}{x_j}
  x_r 
  \,dx_k   
  -
  \int_{\gamma_i}
  \ddt{(s_{jkst}\sigma_{st})}{x_r}
  x_r  
  \,dx_k 
  \Bigr).
  \end{align*}
Thus
  \begin{align*}
  \lrp{ \int_{\gamma_i} \omega \,dx }_j 
  & = l_i\omega_{ji}(\alpha) + 
  \sum_{ \substack{ k,r,s,t\\r\neq j}}\Bigl( 
  \int_{\gamma_i}  
  \ddt{(s_{rkst}\sigma_{st})}{x_j}
  x_r 
  \,dx_k   
  -
  \int_{\gamma_i}
  \ddt{(s_{jkst}\sigma_{st})}{x_r}
  x_r  
  \,dx_k 
  \Bigr).
  \end{align*}
This together with \eqref{eq:auxqp} give the asserted equality.
\end{proof}

To compute the average stress we use the following formulas.

\begin{lemma}\label{lm:averagestresslemma}
Let the local elasticity tensor $c_{ijkl}$ be smooth on $Q$ and $Y$-periodic.
Suppose that $u$ is a quasiperiodic minimizer of  
  \mbox{$\int_Q \e(u)\!\cdot\!\sigma(u) \,dx$.}
Then for $i \neq j$,
  \begin{align*}
  \int_Q \sigma_{ii} \,dx & = l_i \int_{\gamma_j} (\sigma\nu)_i \,dx, &
  \int_Q \sigma_{ij} \,dx & = l_j \int_{\gamma_i} (\sigma\nu)_i \,dx.
  \end{align*}
\end{lemma}
\begin{proof}
Since $u$ is quasiperiodic and minimizes the elastic energy, 
  $u$ is the unique solution, up to translation, to the equation
  \begin{align*}
  \mop{div}\sigma & = 0 \text{ in $\mathcal{D}'(\widetilde{Q})^2$,}\\
  \sigma \nu & = 0 \text{ in $H^{-1/2}(\partial \widetilde{Q})^2$.}
  \end{align*}

Since $\sigma$ and $\mop{div}\sigma$ have components in $L^2(Q')$,
  it follows from the Green formula that
  \begin{align*}
  \int_{Q'} (\sigma_{11}, \sigma_{12}) \cdot \nabla \varphi \,dx
  +
  \int_{Q'} \mop{div} (\sigma_{11}, \sigma_{12}) \, \varphi \,dx
  & = 
  \lran{ (\sigma\nu)_1, \varphi },
  \end{align*}
  for any $\varphi \in H^1(Q')$.
By the regularity of $u$ and the vanishing of $\sigma \nu$ on
  $\partial \widetilde{Q}$, we have that
  $\sigma \nu \in L^2(\partial Q')^2$.
See~\cite{fichera}.
Hence
  \begin{align*}
  \lran{ (\sigma\nu)_1, \varphi } 
  & =  
     \Idx{\gamma_1'}{ (\sigma\nu)_1 \, (\varphi(x) - \varphi(x_1, x_2 - l_2)) }  \\
  & \quad + \Idx{\gamma_2'}{ (\sigma\nu)_1 \, (\varphi(x) - \varphi(x_1 - l_1, x_2)) },
  \end{align*}
  where we in the last step made a change of variables and used the
  periodicity of $\sigma$ and the opposite signs of the outward unit normals.

Let $\varphi = x_1$.
Then $\varphi(x_1, x_2 - l_2) = \varphi(x)$ and
  $\varphi(x) - \varphi(x_1 - l_1,x_2) = l_1$.
Thus
  \begin{align*}
  \lran{ (\sigma\nu)_1, \varphi }
  & = 
    l_1 \Idx{\gamma_2'}{(\sigma\nu)_1}
  = 
    l_1 \Idx{\gamma_2}{(\sigma\nu)_1},
  \end{align*}
  where we in the last step used the periodicity of $\sigma$.
Since $\nabla \varphi = (1,0)$ and $\mop{div}\sigma = 0$, we have by
  periodicity and the Green formula that
  \begin{align}\label{eq:avaux1}
  \Idx{Q}{\sigma_{11}}
  & =
      \Idx{Q'}{\sigma_{11}}
    =
      l_1 \Idx{\gamma_2}{(\sigma\nu)_1}.
  \end{align}

With $\varphi = x_2$,
we obtain in the same way,
  \begin{align}\label{eq:avaux2}
  \Idx{Q}{\sigma_{12}}
  & =
    l_2 \Idx{\gamma_1}{(\sigma\nu)_1}. 
  \end{align}

By the symmetry with respect to the indices, the equations corresponding to \eqref{eq:avaux1} and~\eqref{eq:avaux2} are 
  \begin{align*}
  \Idx{Q}{\sigma_{22}} & = l_2 \Idx{\gamma_1}{(\sigma\nu)_2}, &
  \Idx{Q}{\sigma_{12}} & = l_1 \Idx{\gamma_2}{(\sigma\nu)_2},
  \end{align*}
  which completes the proof.
\end{proof}

\begin{proof}[Proof of Lemma \ref{lm:qpspecial}]
By Lemma~\ref{lm:quasiperiodlemma} we have
\begin{align*}
\xi_{11} l_1 & = \sum_{r,s,t} \int_{\gamma_1} \!\! s_{1rst} \sigma_{st} \,dx_r 
  + \sum_{q,s,t} \lrp{ \int_{\gamma_1} \!\! s_{2qst} \ddt{\sigma_{st}}{x_1}x_2 \,dx_q 
  -  \int_{\gamma_1} s_{1qst} \ddt{\sigma_{st}}{x_2}x_2 \,dx_q }\!.
\end{align*}
By Lemma~\ref{lm:averagestresslemma} we have
  \begin{align*}
  \Inv{l_2}\int_Q \sigma_{22} \,dx & = \int_{\gamma_1} (\sigma\nu)_2 \,dx
      = \int_{\gamma_1} \sigma_{22} \,dx_1 -\int_{\gamma_1} \sigma_{12} \,dx_2, \\
  \Inv{l_1}\int_Q \sigma_{11} \,dx
    & = \int_{\gamma_2} (\sigma\nu)_1 \,dx
      = \int_{\gamma_2} \sigma_{11} \,dx_2 - \int_{\gamma_2} \sigma_{12} \,dx_1.
  \end{align*}
The representation of $\xi_{11}$ is obtained after a short calculation using the isotropy:
  \begin{align*}
  \xi_{11} l_1
    & =  s_{1111}\int_{\gamma_1}  ( (\mop{tr}\sigma,0) - x_2 \ocurl \mop{tr}\sigma ) \,dx
       - \frac{2s_{1212}}{l_2} \int_Q  \sigma_{22} \,dx, 
  \end{align*}
  where we in the last step used $\mop{div}\sigma = 0$.
Similarly, we find
  \begin{align*}
  \xi_{22} l_2 & =  s_{1111}\int_{\gamma_2}  ( (0,\mop{tr}\sigma) + x_1 \ocurl \mop{tr}\sigma ) \,dx
       - \frac{2s_{1212}}{l_1} \int_Q  \sigma_{11} \,dx,
  \end{align*}
  from which the representation of $\xi_{22}$ follows.

We turn to $\xi_{12}$.
Let $\alpha$ denote the common point of $\gamma_1$ and $\gamma_2$.
By Lemma~\ref{lm:averagestresslemma} and the periodicity we have
  \begin{align*}
  \Inv{l_2}\int_Q \sigma_{12} \,dx & = \int_{\gamma_1'} (\sigma\nu)_1 \,dx
      = \int_{\gamma_1'} \sigma_{12} \,dx_1
        - \int_{\gamma_1'} \sigma_{11} \,dx_2, \\
  \Inv{l_1}\int_Q \sigma_{12} \,dx & = \int_{\gamma_2'} (\sigma\nu)_2 \,dx
      = \int_{\gamma_2'} \sigma_{12} \,dx_2
        - \int_{\gamma_2'} \sigma_{22} \,dx_1.
  \end{align*}
By periodicity it follows from Lemma~\ref{lm:quasiperiodlemma} that
  \begin{align*}
  (\xi_{12} + \omega_{21}(\alpha))l_1 
  & =  s_{1111}\int_{\gamma_1'}   ( (0,\mop{tr}\sigma) + x_1 \ocurl \mop{tr}\sigma ) \,dx
            + \frac{2s_{1212}}{l_2} \int_Q  \sigma_{12} \,dx,\\
  (\xi_{21} + \omega_{12}(\alpha))l_2 
  & =  s_{1111}\int_{\gamma_2'}  ( (\mop{tr}\sigma,0) - x_2 \ocurl \mop{tr}\sigma ) \,dx
            + \frac{2s_{1212}}{l_1} \int_Q  \sigma_{12} \,dx.
  \end{align*}
Since $\xi$ is symmetric and $\omega$ is antisymmetric, we have
  \begin{align*}
  2\xi_{12} & = 
      \frac{s_{1111}}{l_1} \int_{\gamma_1'}  ( (\mop{tr}\sigma,0) + x_1\ocurl \mop{tr}\sigma ) \,dx \\
               & \quad    + \frac{s_{1111}}{l_2} \int_{\gamma_2'}  ( (0,\mop{tr}\sigma) - x_2\ocurl \mop{tr}\sigma ) \,dx 
     + \frac{4s_{1212}}{|Y|} \int_Q  \sigma_{12} \,dx,
  \end{align*}
  which completes the proof.
\end{proof}

\begin{proof}[Proof of Lemma \ref{lm:qpspecial2}]
By considering the differences $\Delta \xi_{ij}$ and carrying out the same calculations as in the proof of Lemma~\ref{lm:qpspecial},
the result is obtained by noting that the supposed shift for isotropic $c_{ijkl}$ means, in particular, that $\Delta c_{1111} = 0$.
\end{proof}

The boundary data in Lemma~\ref{lm:michell} is understood to belong to $H^{-1/2}(\Gamma_i)$ and is naturally assumed to satisfy $\sum_i \lran{\sigma \nu, 1}_{\Gamma_i} = 0$ in either of the cases (i) and (ii).
Here $\lran{\cdot,\cdot}_{\Gamma_i}$ denotes the pairing of $H^{1/2}(\Gamma_i)^2$ and its continuous dual,
  where $\Gamma_i$ are the connected components of $\partial \Omega$ with $\Gamma_0$ being the outer part.
See \cite{adams2003sobolev}.

By the connectedness of $\Omega$, 
an element $v \in L^2(\Omega)^2$ satisfies \mbox{$\mop{div}v = 0$} and $\lran{v \nu, 1}_{H^{1/2}(\Gamma_i)} = 0$ if and only if
  there exists a unique $\phi \in H^1(\Omega)/\rr$ such that $v = \mop{curl}\phi$.
Moreover, the following Helmholtz decomposition holds:
\begin{align}\label{eq:helm}
L^2(\Omega)^2 & = \mop{curl} \{ \phi \in H^1(\Omega) : \phi_{|\Gamma_0} = 0, \, \phi_{|\Gamma_i} = \text{constant} \} \oplus \nabla H^1(\Omega)
\end{align}  
where the sum is direct \cite[Chapter~$2$]{girault}.

A necessary and sufficient condition on symmetric $v_{ij} \in L^2(\Omega)$
  that guarantees the existence of a displacement field $u \in H^1(\Omega)^2$ such that $v = \e(u)$
is due to Donati:
  \begin{align}\label{eq:donaticondition}
  \sum_{i,j} \int_\Omega v_{ij} \varphi_{ij} \,dx = 0,
  \end{align}
  for all symmetric $\varphi$ with components in $C^\infty_0(\Omega)$ and such that $\mop{div}\varphi = 0$.
See~\cite{amrouche, donati, ting}.

\begin{proof}[Proof of Lemma~\ref{lm:michell}]
By the Fr\'echet-Riesz argument, there exists a unique
  $u \in H^1(\Omega)^2/\mop{ker}(\e)$ such that
  \begin{align*}
  \odiv \sigma & = 0 \text{ in $\mathcal{D}'(\Omega)^2$,} \\
  \sigma \nu & = g_i \text{ in $H^{-1/2}(\Gamma_i)^2$,}
  \end{align*}
  for any such $g_i$ satisfying $\lran{ g_i, 1 }_{\Gamma_i} = 0$.
Let $\sigma$ be the solution to the above problem for some given $g_i$ with
  constant and isotropic elasticity tensor on $\Omega$.
Then $\sigma$ has components in $L^2(\Omega)$ and thus there exists
  $\phi \in H^1(\Omega)^2$ such that 
  \begin{align}\label{eq:maux}
  \sigma_{11} & = \ddt{\phi_1}{x_2}, &
  \sigma_{22} & = -\ddt{\phi_2}{x_1}, &
  \sigma_{12} & = -\ddt{\phi_1}{x_1} = \ddt{\phi_2}{x_2},
  \end{align}
where the last equation comes from the symmetry of $\sigma$.
The condition~\eqref{eq:donaticondition} can be written as follows in terms of
  $\phi$, for any symmetric $\varphi$ with components in $C_0^\infty(\Omega)$
  satisfying $\odiv \varphi = 0$,
  \begin{align*}
  0 = -\sum_{i,j} \int_{\Omega}  \e_{ij} \varphi_{ij} \,dx 
  & = s_{1111} \int_\Omega \mop{curl}\phi \mop{tr}\varphi \,dx.
  \end{align*}
Indeed,
  \begin{align*}
  \sum_{i,j} \e_{ij}\varphi_{ij}
  & = s_{1111}\lrp{\ddt{\phi_1}{x_2} - \ddt{\phi_2}{x_1}}(\varphi_{11}+\varphi_{22}) \\
  & \quad + 2s_{1212} \lrp{ \ddt{\phi_2}{x_1}\varphi_{11} + \ddt{\phi_2}{x_2}\varphi_{12} - \ddt{\phi_1}{x_1}\varphi_{12}   - \ddt{\phi_1}{x_2}\varphi_{22}},
  \end{align*}
  where we in the last step used \eqref{eq:maux}
  and the isotropy. 
Moreover,
  \begin{align*}
  \int_{\Omega} \lrp{ \ddt{\phi_2}{x_1}\varphi_{11} + \ddt{\phi_2}{x_2}\varphi_{12} } dx
  & = 
  \int_{\Omega} \lrp{ \ddt{\phi_1}{x_1}\varphi_{12} + \ddt{\phi_1}{x_2}\varphi_{22} } dx
  = 0,
  \end{align*}
  where both integrals vanish because $\mop{div}\varphi = 0$.
Thus
  \begin{align*}
  \sum_{i,j} \int_{\Omega} \e_{ij}\varphi_{ij} \,dx = s_{1111} \int_{\Omega} \lrp{ \ddt{\phi_1}{x_2} - \ddt{\phi_2}{x_1} } (\varphi_{11} + \varphi_{22}) \,dx,
  \end{align*}
  as claimed.  
Hence
  $
  \int_\Omega \ocurl \phi \mop{tr} \varphi \,dx = 0,
  $
  which by the sufficiency of the Donati condition~\eqref{eq:donaticondition},
  guarantees the existence of
  some $v \in H^1(\Omega)^2$ such that the given $\sigma$ comes from $v$.
Since this condition does not depend on the elasticity tensor, the
  displacement field $v$ for any supposed elasticity tensor can be recovered
  by using the same $\phi$, which proves the first assertion (i).

We turn to (ii) and note that by isotropy $\Delta s_{1111} = 0$.
By the calculation in (i), 
  \begin{align*}
  \Delta \sum_{i,j} \e_{ij}\varphi_{ij} & = 
  2\Delta s_{1212} \lrp{ \ddt{\phi_2}{x_1}\varphi_{11} + \ddt{\phi_2}{x_2}\varphi_{12} - \ddt{\phi_1}{x_1}\varphi_{12}   - \ddt{\phi_1}{x_2}\varphi_{22}}.
  \end{align*}
Since $\Delta s_{1212}$ is constant, we conclude that
  \begin{align*}
  \Delta \sum_{i,j} \int_{\Omega} \e_{ij}\varphi_{ij} \,dx & = 0,
  \end{align*}
  which shows the claimed independence.

We now consider the case when the condition on the data is relaxed,
  but $c_{ijkl}$ is constant and $\Delta (s_{1122}/s_{1111}) = 0$.
Since \mbox{$\sum_i \lran{\sigma\nu, 1}_{\Gamma_i} = 0$},
  the stress field $\sigma$ exists with components in $L^2(\Omega)$ and by the Helmholtz
  decomposition~\eqref{eq:helm} there exist $\phi,\, q \in H^1(\Omega)^2$ such that
  \begin{align*}
  \mat{ \sigma_{11} & \sigma_{12} \\ \sigma_{12} & \sigma_{22} } 
  =
  \mat{ \ddt{\phi_1}{x_2} & -\ddt{\phi_1}{x_1} \\ \ddt{\phi_2}{x_2} & - \ddt{\phi_2}{x_1}  }
  +
  \nabla q.
  \end{align*}
By the necessity of \eqref{eq:donaticondition}, we have for any $\varphi$ as
  above, that
  \begin{align}\label{eq:donatiaux}
  0
  =
  - \sum_{i,j}\int_\Omega \e_{ij} \varphi_{ij} \,dx
  =
  s_{1111} \int_\Omega \mop{curl}\phi \mop{tr}\varphi \,dx
  -
  s_{1122} \int_\Omega \mop{div}q \, \mop{tr}\varphi \,dx.
  \end{align}
Provided $\lran{\sigma \nu, 1}_{\Gamma_i} \neq 0$ for some $i$, 
  \eqref{eq:donatiaux} holds independently of the elasticity tensor
  if and only if $\Delta (s_{1122}/s_{1111}) = 0$, since $\mop{div}q = 0$ gives self-equilibriated holes, $\lran{\sigma \nu, 1}_{\Gamma_i} = 0$.
Thus, in general, then and only then there exists a displacement field
  $v \in H^1(\Omega)^2$ coming from $\sigma$, by the Donati condition~\eqref{eq:donaticondition}.
\end{proof}

\bibliographystyle{plain}
\bibliography{refs}

\end{document}